\documentclass[11pt]{amsart}
\usepackage{amscd,amssymb,amsmath}
\usepackage[arrow,matrix]{xy}
\usepackage{graphicx}
\newtheorem{thm}[subsection]{Theorem}

\newtheorem{proposition}[subsection]{Proposition}
\newtheorem{cor}[subsection]{Corollary}

\newtheorem{rk}[subsection]{Remark}
\newtheorem{defn}[subsection]{Definition}

\numberwithin{equation}{section} 

\newcommand{\A}{{\mathcal A}}

\newcommand{\bea}{\begin{eqnarray}}
\newcommand{\eea}{\end{eqnarray}}

\newcommand{\R}{\mathbb{R}}

\begin{document}
\title [ Separable   Quadratic Stochastic Operators]
{ Separable  Quadratic Stochastic Operators }

\author {U.\ A.\ Rozikov, S.  Nazir}
 \address{U.\ A.\ Rozikov\\ Institute of mathematics and information technologies,
Tashkent, Uzbekistan.}
\email {rozikovu@yandex.ru}
 \address{S. Nazir\\ The Abdus Salam international center for theoretical physics
Trieste, Italy.}
\email {snazir@ictp.it }

\begin{abstract} We consider quadratic
stochastic operators, which are separable
as a product of two linear operators. Depending on properties of
these linear operators we classify the set of the separable quadratic stochastic operators:
first class of constant operators, second class of linear and third class of
nonlinear (separable) quadratic stochastic operators. Since the properties of operators from the
first and second classes are well-known, we mainly study properties of the operators of
the third class. We describe some Lyapunov functions of the operators and apply them
 to study $\omega$-limit sets of the trajectories generated by the operators. Also we compare our
 results with known results of the theory of quadratic operators
and give some open problems.

\end{abstract}
\maketitle
\section{Introduction} \label{sec:intro} The history of the quadratic
stochastic operators
can be traced back to work of S.Bernshtein\cite{ber}. During more than
85 years this theory developed and many papers were published (see e.g. \cite{ber}-\cite{gar}, \cite{rz}-\cite{rs}). In
recent years it has again become of interest in connection with
numerous applications to many branches of mathematics, biology and
physics.

A quadratic stochastic operator (QSO) has meaning of a population
evolution operator, which arises as follows: Consider a population
consisting of $m$ species. Let $x^0=(x_1^0,...,x_m^0)$ be the
probability distribution of species in the initial generations,
and $P_{ij,k}$ the probability that individuals in the $i$th
and $j$th species interbreed to produce an individual $k$.
Then the probability distribution
$x'=(x_1',...,x'_m)$ of the species in the first
generation can be found by the total probability i.e.

\begin{equation}\label{1}x'_k=\sum_{i,j=1}^{m}P_{ij,k}x_i^0x_j^0,\ k=1,...,m\end{equation}
 where the cubic matrix $P\equiv P(V)=(P_{ij,k})_{i,j,k=1}^m$ satisfying the following conditions

 \begin{equation}\label{2}P_{ij,k}\geq 0,\ \ \sum_{k=1}^m P_{ij,k}=1,\ i,j\in \{1,...,m\}.
\end{equation}

 This means that the association $x^0\rightarrow x'$ defines a map $V$ of the simplex
  \begin{equation}S^{m-1}=\{x=(x_1,...,x_m)\in \R^m\ :\ x_i\geq 0, \sum_{i=1}^m x_i=1\} \end{equation}into itself, called the evolution operator.
Note that each element $x\in S^{m-1}$ is a probability distribution on $E=\{1,...,m\}$.\\ The population evolves by starting from an arbitrary state $x^{(0)}$, then passing to the state $x'=V(x^{(0)})$, then to the state $x''=V(V(x^{(0)})), $and so on.\\ For a given $x^{(0)}\in S^{m-1}$, the trajectory $\{x^{(m)}\ |\ m=0,1,...\}$ induced by the QSO (\ref{1}) is defined by $$x^{(m+1)}=V(x^{(m)}),\textrm{where}\  m=0,1,...$$
One of the main problems in mathematical biology is to study the asymptotic behavior of the trajectories. This problem was solved completely for the Volterra QSO's (see \cite{gan1}, \cite{gan2},\cite{ge}) defined by the relations (\ref{1}),(\ref{2})
, and by the additional assumption \begin{equation}\label{4} P_{ij,k}=0 \ \textrm{for}\ k\in\{i,j\}.\end{equation} The biology meaning of relation (\ref{2}) is obvious: every individual repeats the genotype of one of its parents. For the Volterra QSO the general formula was given in \cite{gan1},
\begin{equation}\label{5}
x_k'=x_k(1+\sum_{i=1}^ma_{ki}x_i),
\end{equation}
where $a_{ki}=2P_{ik,k}-1$ for $i\neq k$ and $a_{kk}=0.$ Moreover, $a_{ki}=-a_{ik}$ and $|a_{ki}|\leq 1.$

In \cite{gan1,gan3}, the theory of QSO (\ref{5}) was developed using theory of the Lyapunov functions and tournaments. But non-Volterra QSOs (which does not satisfy the condition (\ref{4})) were not completely studied. Because there is no general theory that can be applied for study of non-Volterra operators.\\

Each quadratic operator $V$ can be uniquely defined by the cubic matrix $$P\equiv P(V)=\{P_{ij,k}\}_{i,j,k=1}^m$$
satisfying the conditions (\ref{2}). A constructive description of $P$ was given in \cite{gr}, \cite{gan4}. This construction depends on the probability measure $\mu$ which is given on a fixed graph $G$. In \cite{gan4}, it is proved that the QSO resulting from this construction is of Volterra type if and only if the graph $G$ is connected. A construction of QSOs involving a general finite graph and probability measure $\mu$ (here $\mu$ is the product of measures defined on maximal subgraphs of the graph $G$) and yielding a class of non-Volterra operator was described in \cite{rs}. It was shown that if $\mu$ is given as the product of probability measures, then the corresponding non-Volterra operator can be studied by $N$ Volterra operators (where $N$ is the number of connected components of the graph).\\

In this paper we consider QSO (\ref{1}),(\ref{2}) with additional condition
\begin{equation}\label{21}
P_{ij,k}=a_{ik}b_{jk},\ \ \textrm{for all}\ i,j,k\in E
\end{equation}
where $a_{ik},b_{jk}\in \R$ entries of matrices $A=(a_{ik})$ and $B=(b_{jk})$
 such that the conditions (\ref{2}) are satisfied for the coefficients (\ref{21}).

 Then the QSO $V$ corresponding to the coefficients (\ref{21}) has the form
 \begin{equation}\label{22}
x'_k=(V(x))_k=(A(x))_k \cdot(B(x))_k,
\end{equation}
where $(A(x))_k=\sum_{i=1}^ma_{ik}x_i,\ (B(x))_k=\sum_{j=1}^mb_{jk}x_j.$\\
\begin{defn}
The QSO (\ref{22}) is called separable quadratic stochastic operator (SQSO).

\end{defn}
\begin{rk}
{\rm 1.} If $A$ (or $B$) is the identity matrix then the operator (\ref{22}) becomes  a linear  Volterra QSO.

{\rm 2.} The following example shows that the condition (\ref{21}) is sufficient for a QSO to be product
of two linear operators, but the condition is not necessary:
consider matrices
    $$   A=\left( \begin{array}{ccc}
             0 & 1 & 0 \\
             1 & 0 & 0 \\
             0 & 0 & 1 \\
           \end{array}
         \right),\ \ B=\left(\begin{array}{ccc}
             0 & 1 & 0 \\
             1 & 2 & 2 \\
             0 & 2 & 1 \\
           \end{array}
         \right).$$
Then corresponding QSO is
 $$\left\{ \begin{array}{lll}
             x'_1=x^2_2 \\
             x_2'=x_1(x_1+2x_2+2x_3)\\
             x_3'=x_3(2x_2+x_3)\\
           \end{array}\right.$$
For this operator we can take $P_{13,2}=P_{31,2}=1$ but from matrices
$A$ and $B$ we have $a_{32}b_{12}=0\ne P_{31,2}.$
But one can easily check that a QSO with coefficients $P_{ij,k}$
can be written as the product of two linear operators
$A=(a_{ik})$ and $B=(b_{jk})$
if and only if
\begin{equation}\label{o}
P_{ij,k}+P_{ji,k}=a_{ik}b_{jk}+a_{jk}b_{ik}\ \ {\rm for\ any} \ \ i,j,k\in E.
\end{equation}
Thus the condition (\ref{21}) is a particular case of (\ref{o}). In this paper for simplicity
we assume that (\ref{21}) holds.
\end{rk}
The paper is organized as follows. In section 2 we classify the set of the SQSOs as constant operators, linear and nonlinear operators. Since the theory of linear operators is well-known,  we mainly study properties of the nonlinear SQSOs. Section 3 is devoted to describe
some Lyapunov functions of such nonlinear operators. In section 4 we apply the Lyapunov functions to
obtain upper estimates for the set of $\omega$-limit points of trajectories generated by SQSOs. In the last section
we compar our results with known results of QSOs
and give some open problems.
  \section{Classification of SQSO's}

  From the conditions $P_{ij,k}\geq 0$ and  $\sum_{k=1}^m P_{ij,k}=1$ for all $i,j$ it follows
  the condition on matrices $A$ and $B$ that $a_{ik}b_{jk}\geq 0$, $AB^T=\textbf{1}$, where $B^T$ is the transpose of $B$ and  $\textbf{1}$
  is the matrix with all entries $1$'s. If $a^{(i)}=(a_{i1},...,a_{im})$ is the $i$-th row of the matrix $A$ and $b^{(j)}=(b_{j1},...,b_{jm})$ is the $j$-th row of the matrix $B$, then from $AB^T=\textbf{1}$ we get
  \begin{equation}\label{23}
  a^{(i)}b^{(j)}=1,\ \ \ \textrm{for all}\ i,j=1,...,m,
  \end{equation}
  For a fixed $j$, the above condition implies that
  \begin{equation}\label{2a} A(b^{(j)})^T=(1,1,...,1).
  \end{equation}

  If $\det(A)\neq 0$, then (\ref{2a}) gives $b^{(j)}=b^{(k)}$ for all $k,j=1,..,m$, i.e., all rows of $B$ are the same, therefore $\det(B)=0.$ Similarly, if $\det(B)\neq 0$, then  all the rows of $A$ must be the same, so $\det(A)=0.$\\

  But in spite of both $A$ and $B$ have zero determinant we may have the
  matrices $A$ and $B$ are without all identical rows.
  \begin{rk}
  If $m=2 $ and $\det(A)=0$ satisfying the condition (\ref{23}) then one can check that $A$ has the
  identical rows. But for $m\geq 3$, from $\det(A)=0$ and the  condition (\ref{23}),
  we may have $A$ and $B$ having not all the same rows. For example, consider
$$   A=\left( \begin{array}{ccc}
             b & y_1 & 1-y_1 \\
             b & y_2 & 1-y_2 \\
             b & y_3 & 1-y_3 \\
           \end{array}
         \right),\ \ B=\left(\begin{array}{ccc}
             0 & 1 & 1 \\
             0 & 1 & 1 \\
             \frac{1}{2b} & \frac{1}{2} & \frac{1}{2} \\
           \end{array}
         \right),$$ where $b>0$, $y_i\in [0,1]$ for all $i=1,2,3$.
It is easy to see that these matrices satisfy the condition (\ref{23}) and $\det(A)=\det(B)=0$.
  \end{rk}

  So we have three cases for the SQSO (\ref{22}):\\

  \textbf{Case 1:} If $\det(A)=\det(B)=0$ and both of them have the identical rows,
  then the SQSO becomes constant i.e
\begin{equation}\label{c}
x_k'=a_{1k}b_{1k},\ \  \textrm{for }\ \ k=1,...,m.
  \end{equation}

   In this case the dynamical system is trivial: independently on initial
   point $x^{(0)}=x^0\in S^{m-1}$ all trajectories coincide with $\{x^{(n)}\}$ such that $x^{(n)}_k=a_{1k}b_{1k},\ \  \textrm{for }\ \ k=1,...,m, \ \ n=1,2,...$\\

 \textbf{ Case 2:} If $\det(A)\neq 0$, then $B$ has the same rows and SQSO becomes
\begin{equation}\label{l} x_k'=b_{1k}\sum_{i=1}^ma_{ik}x_i,\ \  \textrm{for }\ \ k=1,...,m.\end{equation}
  which is a linear stochastic operator.

  \begin{rk}
{\rm 1.} Since $B$ can be uniquely determined by a given $A$
with $\det(A)\neq 0$, the operator (\ref{l}) depends on $A$ only. Moreover, the matrix $\mathbf{P}=\left(b_{1k}a_{ik}\right)_{i,k=1}^m$ is a quadratic stochastic matrix.

{\rm 2.} It is known (see e.g \cite{sh}) that the properties of homogeneous Markov chains with the phase space
$E=\{1,...,m\}$ can by completely determined by the initial distribution $x\in S^{m-1}$ and
the stochastic matrix $\mathbf{P}$ i.e. the dynamical system generated by the operator (\ref{l}). The theory of
such dynamical systems is known (see for example \cite{sh}). Also for the dynamical behavior of general linear
operators $\mathbf{R}^m\to \mathbf{R}^m$ see, for example, \cite{d}, pages 159-181.
\end{rk}

  \textbf{Case 3:} If $\det(A)=\det(B)=0$ but  both of them don't have all the identical rows, then the SQSO is
  \begin{equation}\label{24}
  x_k'=\left(\sum_{i=1}^ma_{ik}x_i\right)\left(\sum_{j=1}^mb_{jk}x_j\right),\ \  \textrm{for }\ \ k=1,...,m.
  \end{equation}
\begin{rk} {\rm 1.} Since $B$ can {\it not} be uniquely determined by a given $A$
with $\det(A)=0$, the operator (\ref{24}) depends on both matrices $A$ and $B$.

{\rm 2.} By the above mentioned reasons only SQSO (\ref{24}) is interesting to study.
In this case, we can have a rich theory of such operators: to find Lypunov functions;
to study fixed points; to determine concepts of tournaments and so on. So the sequel
of this paper is devoted to SQSO (\ref{24}).  \end{rk}

\section{Lypunov Functions of SQSO (\ref{24})}

Let $x^{(0)}\in S^{m-1}$ be the initial point, and let $\{x^{(0)},x^{(1)},x^{(2)},...\}$ be the trajectory of the point $x^{(0)}$. Denote by $\omega(x^{(0)})$ the set of limit points of the trajectory $\{x^{(n)}\}_{n=0}^\infty$. Since $\{x^{(n)}\}_{n=0}^\infty\subset S^{m-1}$ and $S^{m-1}$ is a compact set, it follows that $\omega(x^{(0)})\neq \emptyset.$ If $\omega(x^{(0)})$ consists of a single point, then the trajectory converges, and $\omega(x^{(0)})$ is a fixed point of the operator $V.$
\begin{defn}
A continuous functional $\phi: S^{m-1}\longrightarrow \R$ is called a Lyapunov function for the dynamical system (\ref{24}) if the limit $\lim_{n\rightarrow \infty}\phi(x^{(n)})$ exists for any initial point $x^{(0)}\in S^{m-1}.$

\end{defn}
Obviously, if $\lim_{n\to\infty}\phi(x^{(n)}) = c$, then $\omega(x^0)\subset \phi^{-1}(c)$. Consequently, for an
upper estimate of $\omega(x^0)$ we should construct a set of Lyapunov functions that is as large
as possible.

Denote

$$\A=\{(A,B)\ :\ \det(A)=\det(B)=0,\ AB^T=\mathbf{1}, \ \textrm{and}$$ $$ \textrm{both matrices don't have all the identical rows}\}.$$
\begin{thm}\label{t1}
For the dynamical system (\ref{24}), the function $\psi_c:S^{m-1}\to \mathbf{R}$ defined by
\begin{equation} \label{31}
\psi_c(x)=\sum_{k=1}^m c_kx_k
\end{equation}
is a Lypunov function if $c=(c_1,...,c_m)^T $ satisfies $c_i\geq 0$ for all $1\leq i \leq m$ and either $Ac\leq Ic$ or $Bc\leq Ic$ where $A=(a_{ij})$, $B=(b_{ij})$, with  $0\leq a_{ij}, b_{ij}\leq 1$ for all $1\leq i,j \leq m$,  $(A,B)\in \A$ and $I$ is the identity matrix of order $m$.
\end{thm}
\proof
Suppose $Ac\leq Ic$, then we have

\begin{eqnarray*}
\psi_c(x')=\sum_{k=1}^mc_kx_k'=\sum_{k=1}^mc_k\sum_{i,j=1}^ma_{ik}b_{jk}x_ix_j\\
\leq \sum_{k=1}^mc_k\sum_{i,j=1}^ma_{ik}x_ix_j=\sum_{k=1}^mc_k\sum_{i=1}^ma_{ik}x_i\\
=\sum_{i=1}^m\left(\sum_{k=1}^mc_ka_{ik}\right)x_i\leq \sum_{i=1}^mc_ix_i=\psi_c(x)
\end{eqnarray*}
as $0\leq b_{ij}\leq 1$.

Thus, for any $n$, we have $\psi_c(x^{(n)})\leq \psi_c(x^{(n-1)})$ and $\underline{c}\leq \psi_c(x^{(n)})\leq \overline{c},$ with $\underline{c}=\min_{i}c_i$, $\overline{c}=\max_ic_i$. Consequently, the sequence $\{\psi_c(x^{(n)})\}_{n=0}^\infty$ is convergent.
Therefore, $\psi_c(x)=\sum_{k=1}^m c_kx_k$ is a Lyapunov function for the dynamical system (\ref{24}).
\endproof

\begin{cor}
The function defined by \begin{equation}\label{34}
\phi(x)=\prod_{k=1}^m\left(\sum_{i=1}^mc_{ik}x_i\right)^{p_k}
\end{equation}
is a Lyapunov function for the dynamical system (\ref{24})
for any $p_k\in \R^+$ if $c^{(k)}=(c_{k1}, c_{k2}, ...,c_{km})$ satisfies $c_{ij}\geq 0$, for all $i,j=1,...,m$ and either $Ac^{(k)}\leq c^{(k)}$ or $Bc^{(k)}\leq c^{(k)}$ for all $k=1,...,m$, where $A=(a_{ij})$
and $B=(b_{ij})$, $(A,B)\in \A$ with $0\leq a_{ij},b_{ij}\leq 1$ for all $1\leq i,j \leq m.$
\end{cor}
\begin{rk}\label{r} To use Theorem \ref{t1} one has to find non-zero solution (if exists) of the system
of inequalities  $Ac\leq Ic$ (or $Bc\leq Ic$). It is known (see \cite{k}, page 42) that
the set $\mathcal{C}$ of all solutions of the system of inequalities is a polyhedral convex cone:
if $c$ is in $\mathcal{C}$, then the vector $tc$ for all $t\geq 0$ ("the ray or halfline generated
by $c$) are also in $\mathcal{C}$. Such a polyhedral convex cone $\{c: (A- I)c\leq 0\}$ can also be expressed
as the convex-cone hull $Q^<$ of a finite set $Q=\{Q_1,...,Q_q\}:$
$$Q^<=\{c: c=v_1Q_1+...+v_qQ_q,\ v_i\geq 0\}=\{c: c=Qv, v\geq 0\}.$$
Concretely,  $Q$ is the $m$ by $q$ matrix with $Q_1$, ...,$Q_q$ as columns.
For a method of solving of the system of linear inequalities see \cite{p}.
\end{rk}
Now we shall give some solutions of the above mentioned system of inequalities.
\begin{proposition}
If \begin{equation}\label{32}a_i=\sum_{k=1}^ma_{ik}\leq 1, \ \ \textrm{for all}\ \  i=1,...,m,
\end{equation} then the system of inequalities $Ac\leq Ic$
has nonzero solutions\\ $c(t)=(tc_1, tc_2,...,tc_m)$ with $c_i=a_i$ for all $t> 0$.
\end{proposition}
\proof By Remark \ref{r} it is enough to prove for $c(1)$ i.e $t=1$:
$$\sum_{k=1}^ma_{ik}c_k=\sum_{k=1}^ma_{ik}a_k\leq \sum_{k=1}^ma_{ik}=a_i=c_i.$$\endproof

The following example shows that the above condition (\ref{32})
is not necessary:

{\bf Example.} Consider
 $$A=\left(
  \begin{array}{ccc}
 1 & 0 & 0 \\[1.5mm]
 \frac{1}{3} & \frac{1}{2} & \frac{1}{4} \\[1.5mm]
\frac{ 2}{3} & \frac{1}{4} & \frac{1}{8} \\
\end{array}
\right),\ \
B=\left(
           \begin{array}{ccc}
             1 & \frac{8-3b_1}{6} & b_1 \\[1.5mm]
             1 & \frac{8-3b_2}{6} & b_2 \\[1.5mm]
             1 & \frac{8-3b_3}{6} & b_3 \\
           \end{array}
         \right),
$$ where $\frac{2}{3}\leq b_j\leq 1$ for all $j=1,2,3$.
It is easy to show that $A$ and $B$ satisfy the condition (\ref{23}).
But the condition (\ref{32}) is not satisfied for $A$.
Moreover, the system of inequalities $Ac\leq Ic$ has many non-zero solutions. More precisely
the complete set of solutions is  $$\mathcal{C}=\left\{c=(c_1,c_2,c_3):\ c_1\geq 0, \ c_2\geq
{11\over 9}c_1, \ {16\over 21}c_1+{2\over 7}c_2\leq c_3\leq 2c_2-{4\over 3}c_1\right\}.$$
It also can be observed that there is no nonzero solution for the system of inequalities
$Bc\leq Ic$.

\section{$\omega$-limit set of SQSO}

For (\ref{c}) we have $\omega(x^0)=\{(a_{11}b_{11},...,a_{1m}b_{1m})\}$, $\forall x^0\in S^{m-1}$.
But for (\ref{l}) the set $\omega(x^0)$ depends on $x^0$ and on the properties of the matrix $A$.
The set $\omega(x^0)$ can contain a single point (ergodic case, see \cite{sh}, page 118), it can be a
finite set (non-ergodic or periodic case, see page 121 of \cite{sh}).

In this section using the Lyapunov functions described in the previous section
we shall give upper estimation of $\omega(x^0)$ for SQSO (\ref{24}).
Denote
$$\mathcal{C}=\{c\in \mathbf{R}^{m}: c_i\geq 0, \ c_1+...+c_m>0, \ Ac\leq Ic \ \ {\rm or} \ \ Bc\leq Ic \}.$$
Then by Theorem \ref{t1} we have that $\psi_c$ is a Lyapunov function for any $c\in \mathcal{C}$. That is for any initial point
$x^0\in S^{m-1}$ we have
\begin{equation}\label{lam}
\lim_{n\to \infty}\psi_c(x^{(n)})=\lambda_c(x^0), \ c\in \mathcal{C}.
\end{equation}
Thus $\omega(x^0)\subset \{x\in S^{m-1}: \psi_c(x)=\lambda_c(x^0)\}$ for any $c\in \mathcal{C}$ which
implies
\begin{equation}\label{w}
\omega(x^0)\subset \bigcap_{c\in \mathcal{C}}\left\{x\in S^{m-1}: \psi_c(x)=\lambda_c(x^0)\right\}.
\end{equation}
The estimation (\ref{w}) is very useful: assume that there are $m$ distinct vectors $c^{(1)},...,c^{(m)}\in \mathcal{C}$
such that $\det(C)\ne 0$ where $C$ is the $m\times m$ matrix with rows $c^{(i)}$, $i=1,...,m$. Then
system of equations $\psi_{c^{(i)}}(x)=\lambda_{c^{(i)}}(x^0)$, $i=1,...,m$ has unique solution $x=x^*$
which by (\ref{w}) gives that $\omega(x^0)=\{x^*\}$.  If there is no any collection $c^{(i)}$, $i=1,...,m$ with
$\det(C)\ne 0$, then RHS of (\ref{w}) is a uncountable set. Note that RHS of (\ref{w}) can not be empty set
since $\omega(x^0)\ne \emptyset$, because $\{x^{(n)}\}_{n=0}^\infty\subset S^{m-1}$ and $S^{m-1}$ is a compact set.\\
\section{Discussions} As it was mentioned in the Introduction, $P_{ij,k}$ is the probability
that individuals in the $i$th
and $j$th species interbreed to produce an individual $k$. Our assumption (\ref{21}) means that the parents $ij$ independently participate for producing $k$. Under this condition the QSO became a product of two
linear operators. Note that Volterra QSO (\cite{gan1}) also is product of two linear operators: identical operator and a linear operator. The theory of the Volterra QSOs is developed enough. But our SQSO corresponding to quadratic matrices $A$ and $B$ coincides with Volterra QSO
iff $A$ or $B$ is the identical matrix. Thus one can expect that SQSOs may have different dynamical behavior.
We already have seen that SQSOs (for example, (\ref{l})) may have periodic trajectories this is quite
different behavior from the behavior of the Volterra QSO, since Volterra operators have no periodic
trajectories. We described a wide family of linear Lyapunov functions for SQSOs, but for Volterra QSOs, such a family of functions is not described, except $\varphi(x)=\sum_{i=1}^rx_i$ which is constructed in \cite{gan1} under
some conditions on parameters of (\ref{5}), by using theory of tournaments. We think our argument also works for Volterra QSOs to describe
linear Lyapunov functions.

We know that $\lambda_c(x^0)$ given in (\ref{lam}) exists, but, in general, we don't know how to 
compute them.  The following interesting problems are also open: Is there any other kind of Lyapunov functions of SQSOs?
Develop a tournament theory approach for SQSOs (\ref{24}). Here one may define two tournaments
corresponding to quadratic matrices $A$ and $B$. Then properties of (\ref{24}) will depend
on both tournaments. Also very interesting problem is to study connection between these tournaments.\\

{\bf Acknowledgements.} This work was done in the Abdus Salam International Center for
Theoretical Physics (ICTP), Trieste, Italy. UAR thanks the ICTP  for providing financial
support of his visit (within the scheme of Junior Associate) to ICTP (February-April 2009).
SN thanks the ICTP for providing to her a PostDoc position (January-December 2009).

{}

\begin{thebibliography}{}

\bibitem{ber} S. N. Bernstein, "The solution of a mathematical problem related to the theory of heredity", Uchen. Zapiski Nauchno-Issled. Kafedry Ukr. Otd. Matem. 1, 83-115 (1924).

\bibitem{d} R.L. Devaney, "An introduction to chaotic dynamical system",  Westview Press, 2003.

\bibitem{gr} N. N. Ganikhodjaev and U. A. Rozikov, "On quadratic stochastic operators generated by Gibbs distributions", Regul. Chaotic Dyn. 11(4), 467-473 (2006).

\bibitem{gan4} N. N. Ganikhodjaev, "On the application of the theory of Gibbs distributions in mathematical genetics", Russian Acad. Sci. Dokl. Math. 61(3), 321-323 (2000).

\bibitem{gan1} R. N. Ganikhodzhaev, "Quadratic stochastic operators, Lyapunov functions, and tournaments", Russian Acad. Sci. Sb. Math. 76(2), 489-506 (1993).

\bibitem{gan2} R. N. Ganikhodzhaev, "On the definition of quadratic bistochastic operators", Russian Math. Surveys 48(4), 244-246 (1993).

\bibitem{ge} R. N. Ganikhodzhaev and D. B. Eshmamatova, "Quadratic automorphisms of a simplex and the asymptotic behaviour of their trajectories", Vladikavkaz. Math. Zh. 8(2), 12-28 (2006).

\bibitem{gan3} R. N. Ganikhodzhaev, "Map of fixed points and Lyapunov functions for a class of discrete dynamical systems",  Math. Notes. 56(5), 1125-1131 (1994).

\bibitem{gar} R. N. Ganikhodzhaev and U.A. Rozikov, "Quadratic stochastic operators: results and open problems" arXiv:0902.4207v2 [math.DS]

\bibitem{k} H.W. Kuhn and A.W. Tucker,(eds) "Linear inequalities and related systems", Annal. Math. Stud.
Princeton Univ.Press. 1985.

\bibitem{p} L-P. Pang, E.Spedicato, Z-Q. Xia and W. Wang,"A method for solving the system of linear equations and linear inequalities", Math. Comp. Model. 46, 823-836 (2007).

\bibitem{rz}  U.A. Rozikov and U.U. Zhamilov, "On $F$-quadratic stochastic operators".
 Math. Notes. 83(4), 554-559 (2008).

\bibitem{roz}  U.A. Rozikov and  A. Zada, "On $\ell$-Volterra quadratic stochastic
operators", Doklady Math. 79(1), 32-34 (2009).

\bibitem{rs} U. A. Rozikov and N. B. Shamsiddinov, "On non-Volterra quadratic stochastic operators generated by a     product measure", Stoch, Anal. Appl. 27(2) 353-362(2009).

\bibitem{sh} A.N. Shiryaev, "Probability", 2nd Ed., Springer, 1996.
\end{thebibliography}
\end{document}